\documentclass[14pt]{extarticle}

\usepackage[T2A]{fontenc}
\usepackage[cp1251]{inputenc}
\usepackage[english,russian]{babel}

\usepackage{indentfirst}

\usepackage[tbtags]{amsmath}
\usepackage{amsfonts,amssymb}

\usepackage{mathrsfs}
\usepackage{graphicx}

\numberwithin{equation}{subsection}
\def\hS{\ensuremath{\hat{S}}}
\def\hT{\ensuremath{\hat{T}}}

\def\hR{\ensuremath{\hat{R}}}

\def\Z{\ensuremath{\mathbb{Z}} }
\def\1{{\bf 1} }
\def\2{{\bf 2} }
\def\N{\ensuremath{\mathbb{N}} }

\def\T{\ensuremath{\mathbb{T}^2} }
\def\A{\ensuremath{\mathcal{A}} }

\newcommand{\const}{\mathrm{const}}
\newtheorem{lemma}{Лемма }
\newtheorem{statement}{Утверждение}
\newtheorem{theorem}{Теорема}

\def\proof{\par\noindent{\it Доказательство.}}
\def\endproof{\par\noindent\hbox to \textwidth{\hfill$\blacksquare$}}

\begin{document}

\author{Р.\,А.~Конев,  В.\,В.~Рыжиков\footnote { Работа выполнена в рамках программы  гранта Президента РФ для поддержки ведущих научных школ НШ-5998.2012.1}}

\title
{\large О  наборе спектральных кратностей $\{2, 4, \ldots,2^n\}$\\  для вполне эргодических $\Z^2$-действий}

\maketitle

\begin{abstract}
Работа посвящена  реализации наборов спектральных кратностей для эргодических $\Z^2$-действий.  Даны достаточные условия, обеспечивающие набор кратностей вида $\{2,4,\ldots,2^n\}$.  Использумый метод
позволяет с помощью специальной предельной процедуры  получить аналогичный результат для перемешивающих действий. 

Библиография: 12 названий.
\end{abstract}


\subsection{Введение}
\label{subsec1}

Задача о реализации заданного подмножества натурального ряда как множества кратностей спектра эргодического действия хорошо известна  в спектральной теории динамических систем с инвариантной мерой (см. \cite{Anosov},\cite{Goodson}). Современное состояние этой тематики и методология полно отражены  в недавнем обзоре \cite{Danilenko_survey_2011}.
Для  $\Z$-действий в  \cite{Ryzhikov}, \cite{Ageev} независимо был получен простейший набор $\{2\}$, что было ответом на вопрос Рохлина о существовании эргодического автоморфизма с  однородным спектром кратности два. Авторы настоящей статьи позднее   нашли реализацию  наборов кратностей $\{2, 4,\ldots,2^n\}$ для $\Z$-действий (см. \cite{Ryzhikov_MatSbornik_2009}).  Общие результаты были предложены  Катком и Леманчиком \cite{Katok},   Даниленко и Соломко (см. \cite{Danilenko_survey_2011},\cite{Solomko}).

Настоящая  работа посвящена реализации спектральных наборов кратности $\{2,4,\ldots,2^n\}$ для эргодических $\Z^2$-действий. Используя методы, предложенные  в \cite{Ryzhikov_WeekLimits}, этот результат можно применить для построения перемешивающего $\Z^2$-действия с такими же наборами спектральных кратностей.      Реалицация однородного спектра кратности $\{2\}$ получена первым из авторов,  остальные результаты  совместные.   Перемешивающие $\Z^2$-действия планируется рассмотреть в отдельной  работе.

В статье   использованы следующие обозначения:
\begin{description}
  \item[$\mathbb{T}$] -- единичная окружность,
  \item[$I$] -- тождественный оператор,
  \item[$\Theta$] -- оператор ортопроекции на подпространство констант,
  \item[$\{\const\}$] -- подпространство констант в $L_2(Y,\nu)$,
  \item[$L_2^0(Y,\nu)$] -- подпространство функций $L_2(Y,\nu)\ominus\{\const\}$.
\end{description}

\vskip 1em

Множество автоморфизмов $\{S^z\}_{z\in\Z^2}$ пространства Лебега
$(X,\A,\mu)$, $\mu(X)=1$, называется \emph{$\Z^2$-действием},
если эти автоморфизмы удовлетворяют свойству:
$$
    S^{z_1}\circ S^{z_2}=S^{z_1+z_2} \text{ для любых } z_1,z_2\in\Z^2.
$$

Каждому автоморфизму $S$ отвечает оператор
$\hS\colon L_2(X,\mu)\to L_2(X,\mu)$, заданный формулой
$$
    (\hS f)(x)=f(Sx).
$$

Подпространство $C\subseteq L_2(X,\nu)$ для действия $\{S^z\}_{z\in\Z^2}$
является \emph{циклическим подпространством}, если найдётся вектор $h\in L_2(Y,\nu)$ такой, что $\overline{\langle \hS^zh \rangle}_{z\in\Z^2}=C$.
Вектор $h$ называют \emph{циклическим вектором}.

Эргодическое действие $\{S^z\}_{z\in\Z^2}$ имеет \emph{простой спектр},
если $L_2^0(Y,\nu)$ является циклическим подпространством для этого действия.

Действие $\{S^z\}_{z\in\Z^2}$ имеет \emph{однородный спектр кратности $n$},
если пространство $L_2^0(Y,\nu)$ можно представить в виде ортогональной суммы
$\bigoplus_{i=1}^{n} C_i$ циклических подпространств $C_i$, причём ограничения
нашего действия на $C_i$ изоморфны между собой.
Число $n$  принимает значения $1,2,\ldots$  и  $\infty$.

В статье даны достаточные условия для $\Z^2$-действий, обеспечивающие набор кратностей $\{2,4,\ldots, 2^n\}$. В первой части статьи мы изучаем $\Z^2$-действия, имеющих однородный спектр кратности два. Во второй части статьи приведены достаточные условия для $\Z^2$-действий, обеспечивающие  набор кратностей $\{2,4\}$. В третьей части статьи доказан  основной результат 
о наборе кратностей $\{2,4,\ldots,2^n\}$.

Нам следует сказать  о примерах, удовлетворяющих достаточным условиям.  Хорошо известно (это явно не опубликовано, но вытекает, например, из результатов работы \cite{Tikhonov}),  что свойства, которые накладываются на каждый сомножитель
в  произведении $\{\hT_1^z\otimes\ldots\otimes T_k^z\}_{z\in\Z^2}$, являются типичными. Отсюда следует существование нужного набора.     Авторы планируют изложить   в отдельной работе построение частично перемешивающих конструкций,
удовлетворяющих достаточным условиям.  Отметим, что при помощи  предельной процедуры, аналогичной той, что использовалась для
$\Z$-действий в \cite{Ryzhikov_WeekLimits}  (см. также \cite{T}),   получаются  перемешивающие $\Z^2$-действия с набором спектральных кратностей $\{2,4,\ldots,2^n\}$.  Непосредственно к самим перемешивающим системам нельзя применить
метод слабых пределов (в силу их тривиальности последних). Но граница изучаемого нами класса действий содержит перемешивающие
системы, которые наследуют нужные спектральные свойства.

\subsection{Однородный спектр кратности два для эргодического $\Z^2$-действия}
\label{subsec2}

Эргодические $\Z$ действие с однородным спектром кратности два были получены  в \cite{Ageev}, \cite{Ryzhikov}.  В дальнейшем мы сосредоточим наш интерес на вполне эргодических $\Z^2$-действиях: каждый
элемент, кроме тождественного, является эргодическим.   Отметим, что в нашей ситуации
каждый такой элемент будет обладать свойством 
слабого перемешивания (обладать непрерывным спектром).

\begin{theorem}\label{theorem_n1}
    Для того, чтобы вполне эргодическое действие $\{T^z\times T^z\}_{z\in\Z^2}$
    имело однородный спектр кратности два, достаточно выполнение следующих условий:
    \begin{enumerate}
        \item действие $\{T^z\}_{z\in\Z^2}$ имеет простой спектр;
        \item действие $\{T^z\}_{z\in\Z^2}$ имеет слабые пределы 
            $$
                a_1I+b_1\hT^{(1,0)}+c_1\Theta, \quad
                a_2I+b_2\hT^{(0,1)}+c_2\Theta, \quad
                a_3I+b_3\hT^{(1,1)}+c_3\Theta,
            $$
            для некоторых наборов  положительных чисел $a_j, b_j, c_j$, $ a_j + b_j + c_j = 1$,  $j=\{1,2,3\}$.
            
    \end{enumerate}
\end{theorem}

Для доказательства  теоремы нам понадобятся следующие вспомогательные  утверждения.
\par\noindent
Пусть $f$ --- циклический вектор для действия $\{S^z\}_{z\in\Z^2}$.
\emph{Спектральная теорема для $\Z^2$-действий} утверждает существование
изометрии $\Phi\colon L^0_2(Y, \nu)\to L_2(\T, \sigma)$ такой, что
$$
    \Phi (\hS^z) \Phi^{-1}g(u)=u^zg(u) = u_1^{z_1}u_2^{z_2}g(u),
$$
для любых $z=(z_1,z_2)\in\Z^2$ и $u=(u_1,u_2)\in\T$. Поэтому для действия
$\{S^z\otimes S^z\}_{z\in\Z^2}$ найдётся изометрия
$\Psi\colon L^0_2(Y\times Y, \nu\times\nu)\to L_2(\T\times \T, \sigma\times\sigma)$
такая, что
$$
    \Psi (\hS^z\otimes\hS^z) \Psi^{-1} g(u,w) = u^zw^zg(u,w) = u_1^{z_1}u_2^{z_2}w_1^{z_1}w_2^{z_2}g(u,w),
$$
для любых $z=(z_1,z_2)\in\Z^2$, $u=(u_1,u_2)\in\T$ и $w=(w_1,w_2)\in\T$.

Далее наше действие будет рассмотрено в спектральном представлении.

Пространство $L_2(\T\otimes\T,\sigma\otimes\sigma)$ представим в виде прямой суммы
двух подпространств $H_s\oplus H_a$, где
\begin{gather*}
    H_s = \{g(u,w)\colon g(u,w)=g(w,u)\}, \\
    H_a = \{g(u,w)\colon g(u,w)=-g(w,u)\}.
\end{gather*}

\begin{lemma} \label{lemma1}
    Действие $\{\hS^z\otimes \hS^z\}_{z\in\Z^2}$ в пространстве $L_2(\T\otimes\T,\sigma\otimes\sigma)$
    имеет спектр кратности два, если выполнены следующие условия:
    \begin{enumerate}
        \item действие $\{S^z\}_{z\in\Z^2}$ не имеет в спектре дискретной компоненты (спектр непрерывен).
        \item ограничение действия $\{\hS^z\otimes \hS^z\}_{z\in\Z^2}$ на подпространство $H_s$ имеет простой спектр.
    \end{enumerate}
\end{lemma}

\proof{} Обозначим ограничение действия $\{\hS^z\otimes \hS^z\}_{z\in\Z^2}$ на подпространство $H_s$ через
$\{\hS^z_s\otimes \hS^z_s\}_{z\in\Z^2}$, а на подпространство $H_a$ через $\{\hS^z_a\otimes \hS^z_a\}_{z\in\Z^2}$.

Покажем, что подпространства $H_s, H_a$
инвариантны относительно действия $\{\hS^z\otimes \hS^z\}_{z\in\Z^2}$,
причём существует обратимый оператор $F\colon H_s\to H_a$ такой, что
$$
    F (\hS^z_s\otimes\hS^z_s) = (\hS^z_a\otimes\hS^z_a) F.
$$

Инвариантность подпространств $H_s$ и $H_a$ вытекает из их определения.
Также отметим, что  $H_s$ и $H_a$ ортогональны.

Для построения оператора $F$ нам понадобится следующие разбиение множества $\T\times\T$.

Обозначим диагональ $\{(u,u)\colon u\in\T\}$ через $\Delta$.
Множество $\T\times\T\setminus\Delta$ представим в виде $Q\sqcup RQ$,
где $R$ ---  преобразование симметрии относительно диагонали $\Delta$, а
$Q$ --- некоторое борелевское множество.

Замечая, что 
$\sigma\otimes\sigma(\Delta)=0$ в силу непрерывности меры $\sigma$ (отсутствует дискретная компонента),
определим оператор $F\colon H_s\to H_a$:
$$
    Fg=(\chi_\Delta + \chi_Q-\chi_{RQ})g = (\chi_Q-\chi_{RQ})g,
$$
где функции $\chi_\Delta$, $\chi_Q$ и $\chi_{SQ}$ является индикаторами соответствующих множеств $\Delta$, $Q$ и $RQ$.
Оператор $F$ является обратимым, так как $F^2g=(\chi_Q - \chi_{SQ})^2g=g$, т.е. $F^{-1} = F$.

Покажем теперь, что $F (\hS^z_s\otimes\hS^z_s) = (\hS^z_a\otimes\hS^z_a) F$. Действительно, пусть
для любой функции $g_s\in H_s$ имеем
\begin{align*}
    F (\hS^z_s\otimes\hS^z_s) g_s(u,w)& = F u^zw^zg_s(u,w) = (\chi_Q-\chi_{RQ}) u^zw^zg_s(u,w) =\\
        & = u^zw^z (\chi_Q-\chi_{RQ})g_s(u,w) = u^zw^z Fg_s(u,w) = \\
        & = (\hS^z_a\otimes\hS^z_a) F g_s(u,w).
\end{align*}

Таким образом,  действие  $\{\hS^z\otimes \hS^z\}$ является прямой суммой двух унитарно изоморфных действий
с простым спектром. Это означает, что оно имеет однородный спектр кратности 2.
\endproof

 Пусть $f\odot g$ обозначает функцию $\dfrac{1}{2} \left[{f(u)g(w) + f(w)g(u)}\right]$.
\begin{lemma}\label{lemma2}
    Линейная оболочка  $\langle f\odot g \colon f,g\in L_2(\T,\sigma)\rangle$
    всюду плотна в пространстве $H_s$.
\end{lemma}

Лемма очевидным образом вытекает из того, что произвольную функцию $g(u,w)\in H_s$ можно приблизить
конечными линейными комбинациями произведений  вида $f(u)h(w)$,
где $f,h\in L_2(\T,\sigma)$.

\begin{statement}\label{statement1}
    Пусть действие $\{S^z\}_{z\in\Z^2}$ обладает слабым пределом $aI + b\hS^m$  для некоторых $a,b>0$.
    Тогда замыканию линейной оболочки 
    $\langle \hS^zf\otimes \hS^zf\;|\; z\in\Z^2 \rangle$ принадлежат все векторы вида
    $\hS^{z+tm}f\odot \hS^zf$, для всех $z\in\Z^2$, $t\in\Z$.
\end{statement}

\proof{} Обозначим через $C$ замыкание линейной оболочки $\langle \hS^zf\otimes \hS^zf\;|\; z\in\Z^2 \rangle$. Пусть $\{m_i\}\colon\hS^{m_i} \to aI + b\hS^m$. Тогда из определения пространства $C$ следует, что векторы вида $\hS^{z\pm
m_i}f\odot \hS^{z\pm m_i}f$ принадлежат $C$, для всех $z\in\Z^2$.
Устремив $i$ к бесконечности, получим, что
\begin{align*}
    \hS^{m_i}&f \otimes \hS^{m_i}f \xrightarrow{i\to\infty}
            (af+b\hS^{ m}f)\otimes (af+b\hS^{m}f) = \\
        & = a^2 f\otimes f + 2ab\hS^{m}f\odot f + b^2\hS^{ m}f\otimes \hS^{ m}f.
\end{align*}
Отсюда $\hS^{ m}f\odot f\in C$.

Рассмотрим векторы вида $\hS^{ m_i + m}f\odot \hS^{ m_i}f\in C$ и перейдем к пределу. 
Получим, что вектор  

$$(a\hS^{ m}f+b\hS^{ 2m}f)\odot (af+b\hS^{m}f) 
        = ab\hS^{ 2 m}f\otimes f  +$$
$$ +a^2\hS^{ m}f\otimes f +  ab\hS^{m}f\odot \hS^{ m}f   + b^2\hS^{ 2m}f\otimes \hS^{ m}f$$
лежит в $C$.   Так как три последних слагаемых    лежат в $C$,  то это верно и для $\hS^{ 2 m}f\otimes f$.

Аналогично, шаг за шагом,  получим
$$
    \hS^{ tm}f\odot f\in C \text{ для всех } t\in\Z_{+}.
$$
Заметим, что из $\hS^{m_i}\xrightarrow{i\to\infty}aI+b\hS^m$ вытекает
$$
    \hS^{-m_i}\xrightarrow{i\to\infty} (aI+b\hS^m)^* = aI + b\hS^{-m}.
$$

Применяя аналогичным образом последовательность $\{-m_i\}$, окончательно
получим
$$
 \hS^{ z+tm}f\odot \hS^zf\in C \text{ для всех } t\in\Z.
$$

\endproof

\begin{lemma}\label{lemma3}
    Действие $\{S^z\times S^z\}_{z\in\Z^2}$ в пространстве $H_s$  имеет простой спектр,
    если выполнены условия:
    \begin{enumerate}
        \item действие $\{S^z\}_{z\in\Z^2}$ имеет простой спектр;
        \item   для некоторых $a_j,b_j > 0$, $j\in\{1,2,3\}$,  действие $\{S^z\}_{z\in\Z^2}$ имеет слабые пределы вида
            $$
                a_1I+b_1\hS^{(1,0)}, \quad a_2I+b_2\hS^{(0,1)}
                \quad\text{ и }\quad a_3I+b_3\hS^{(1,1)}.
            $$
          
    \end{enumerate}
\end{lemma}

\proof{} Пусть для последовательностей $\{w_i\}$, $\{u_i\}$, $\{v_i\}$ выполнено 
$$
    \hS^{w_i}\to a_1I+b_1\hS^{(1,0)}, \quad \hS^{u_i}\to a_2I+b_2\hS^{(0,1)}
    \quad\text{ и }\quad \hS^{v_i}\to a_3I+b_3\hS^{(1,1)}.
$$
Замыкание линейной оболочки
$C=\overline{\langle \hS^zf\otimes \hS^zf\colon z\in\Z^2\rangle}$ содержится
в пространстве $H_s$, где $f$ --- циклический вектор для
действия $\{\hS\}_{z\in\Z^2}$. Если мы покажем, что пространство $C$ содержит
векторы вида $\hS^{z+(s,t)}f\odot \hS^zf$ для всех $z\in\Z^2$ и $s,t\in\Z$, то
с учётом леммы~2 получим, что $C=H_s$.

Применяя утверждение~\ref{statement1} для последовательностей $w_i$, $u_i$ и $v_i$, получим, что векторы
$\hS^{z+(t,0)}f\odot \hS^zf$, $\hS^{z+(0,t)}f\odot \hS^zf$, $\hS^{z+(t,t)}f\odot \hS^zf$ принадлежат пространству $C$
для любого $z\in\Z^2$ и $t\in\Z$.

Так как векторы $\hS^{z+(t,t)}f\odot \hS^{z}f$ для любого $z\in\Z^2$ и $t\in\Z$
принадлежат замыканию линейной оболочки $C$, то
\begin{align*}
    \hS^{w_i+z+(t,t)}f \odot \hS^{w_i+z} &\xrightarrow{i\to\infty}
            \hS^{z+(t,t)}(a_1f+b_1\hS^{(1,0)}f)\odot \hS^z(a_1f+b_1\hS^{(1,0)}f) =\\
        & = a_1^2\hS^{z+(t,t)}f\odot \hS^zf + a_1b_1\hS^{z+(t,t)}f\odot \hS^{z+(1,0)}f + \\
        & + a_1b_1\hS^{z+(t+1,t)}f\odot \hS^zf + b_1^2\hS^{z+(t+1,t)}f\odot \hS^{z+(1,0)}f, \\
    \hS^{u_i+z+(t,t)}f \odot \hS^{u_i+z} &\xrightarrow{i\to\infty}
            \hS^{z+(t,t)}(a_2f+b_2\hS^{(0,1)}f)\odot \hS^z(a_2f+b_2\hS^{(0,1)}f) =\\
        & = a_2^2\hS^{z+(t,t)}f\odot \hS^zf + a_2b_2\hS^{z+(t,t)}f\odot \hS^{z+(0,1)}f + \\
        & + a_2b_2\hS^{z+(t,t+1)}f\odot \hS^zf + b_2^2\hS^{z+(t,t+1)}f\odot \hS^{z+(0,1)}f.
\end{align*}

Получим для любых $z\in\Z^2$ и $t\in\Z$, что
\begin{gather}
    \hS^{z+(t,t)}f\odot \hS^{z+(1,0)}f + \hS^{z+(t+1,t)}f\odot \hS^zf \in C \label{eq1}\text{ и }\\
    \hS^{z+(t,t)}f\odot \hS^{z+(0,1)}f + \hS^{z+(t,t+1)}f\odot \hS^zf \in C \label{eq2}.
\end{gather}

Подставляя $z:=z+(0,1)$ и $t:=t-1$ в выражение~(\ref{eq1}),  получим
\begin{equation}
    \hS^{z+(t-1,t)}f\odot \hS^{z+(1,1)}f + \hS^{z+(t,t)}f\odot \hS^{z+(0,1)}f\in C \label{eq3}
\end{equation}

В силу (\ref{eq2}), (\ref{eq3}) имеем
$$
    \hS^{z+(t-1,t)}f\odot \hS^{z+(1,1)}f - \hS^{z+(t,t+1)}f\odot \hS^zf\in C.
$$
Подставляя $t=0,\pm1,\pm2,\ldots$, получим $\hS^{z+(t-1,t)}f\odot \hS^zf\in C$.

Итак, мы показали, что для любого $z\in\Z^2$ и $t\in\Z$
$$
    \hS^{z+(t-1,t)}f\odot \hS^{z}f \text{ и } \hS^{z+(t,t)}f\odot \hS^{z}f\in C.
$$
Докажем, что векторы вида $\hS^{z+(t\pm i,t)}f\odot \hS^{z}f$ для
любого $i\in\N$ принадлежат $C$. Так как
$\hS^{z+(t-1,t)}f\odot \hS^{z}f\in C$ и $\hS^{z+(t,t)}f\odot \hS^{z}f\in C$,
то
\begin{align*}
    \hS^{w_i+z+(t,t)}f \odot & \hS^{w_i+z}f \xrightarrow{i\to\infty} \\
        & a_1^2\hS^{z+(t,t)}f\odot \hS^zf + a_1b_1\hS^{z+(t,t)}f\odot \hS^{z+(1,0)}f + \\
        & + a_1b_1\hS^{z+(t+1,t)}f\odot \hS^zf + b_1^2\hS^{z+(t+1,t)}f\odot \hS^{z+(1,0)}f, \\
    \hS^{w_i+z+(t-1,t)}f \odot & \hS^{w_i+z}f \xrightarrow{i\to\infty} \\
        & a_1^2\hS^{z+(t-1,t)}f\odot \hS^zf + a_1b_1\hS^{z+(t-1,t)}f\odot \hS^{z+(1,0)}f + \\
        & + a_1b_1\hS^{z+(t,t)}f\odot \hS^zf + b_1^2\hS^{z+(t,t)}f\odot \hS^{z+(1,0)}f.
\end{align*}
Отсюда $\hS^{z+(t+1,t)}f\odot \hS^zf, \hS^{z+(t-2,t)}f\odot \hS^zf \in C$.
Повторяя эти рассуждения к полученным векторам, мы получим, что
$\hS^{z+(t\pm i,t)}f\odot \hS^zf\in C$ для любого $z\in\Z^2$, $t\in\Z$ и $i\in\N$,
что равносильно
$$
    \hS^{z+(s,t)}f\odot \hS^zf\in C \text{ для любого } z\in\Z^2 \text{ и } s,t\in\Z.
$$
Лемма 3 доказана.

Доказательство теоремы 1.   Пусть для последовательностей 
 $\{p_i\}$, $\{q_i\}$, $\{r_i\}$  при $i\to\infty$ выполнено
\begin{gather*}
    \hT^{p_i}\to a_1I+b_1\hT^{(1,0)}+c_1\Theta,\quad
    \hT^{q_i}\to a_2I+b_2\hT^{(0,1)}+c_2\Theta, \\
    \hT^{r_i}\to a_3I+b_3\hT^{(1,1)}+c_3\Theta.
\end{gather*}
Обозначим через $L_2^{x0}(X\otimes X)$ подпространство функций
$g(x, y)\in L_2(X\otimes X,\mu\otimes\mu)$, зависящих только от одной переменной $x$,
т.е. $g(x,y)=\phi(x)\in L_2^0(X,\mu)$. Аналогичным образом определяется подпространство
$L_2^{0y}(X\otimes X)$. Подпространство $L_2^{00}(X\otimes X)$
определяется как пространство функций $g(x,y)\in L_2(X\otimes X,\mu\otimes\mu)$, имеющих
нулевое средние по переменным $x$, $y$, т.е $\int_X g_{x_0}(y)\,d\mu=0$
и $\int_X g_{y_0}(x)\,d\mu=0$, где $g_{x_0}(y)\colon y\to g(x_0,y)$,
и $g_{y_0}(x)\colon x\to g(x,y_0)$.

Заметим, что имеет место ортогональное разложение
$$
    L_2^0(X\otimes X, \mu\otimes\mu) =
    L_2^{x0}(X\otimes X)\oplus L_2^{0y}(X\otimes X)\oplus L_2^{00}(X\otimes X).
$$

Каждое из этих подпространств инвариантно относительно действия $\{\hT^z\otimes \hT^z\}_{z\in\Z^2}$.

Ограничение действия $\{\hT^z\otimes\hT^z\}_{z\in\Z^2}$ на сумму  $L_2^{x0}(X\otimes X)$
и $L_2^{0y}(X\otimes X)$ изоморфно сумме двух действий  $\{\hT^z\}_{z\in\Z^2}$ в пространстве $L_2^0(X,\mu)$.
Это  ограничение имеет однородный спектр кратности 2.  Непрерывность спектра очевидным образом вытекает из
наличия слабого предела вида 
$$
    \hT^{p_i}\to aI+b\hT^{z}+c\Theta
$$
при $c >0$.

Ограничение действия $\{\hT^z\otimes\hT^z\}_{z\in\Z^2}$ на $L_2^{00}(X\otimes X)$
 изоморфно прямой сумме двух копий действия  $\{\hT^z\odot \hT^z\}_{z\in\Z^2}$.
Поэтому тут мы тоже получаем однородную кратность 2.

Осталось  показать, что действия $\{\hT^z\}$  и $\{\hT^z\odot \hT^z\}$ как унитарные представления группы $\Z^2$  дизъюнктны, что равносильно 
отсутствию ненулевого оператора, сплетающего эти действия.
Пусть для всех $z\in Z^2$
$$J\hT^z = (\hT^z\otimes\hT^z)J.$$

Тогда при $m=(1,0)$ (аналогично при $m=(0,1)$) имеем
$$J(aI+b\hT^{m}) = [(aI+b\hT^{m})\otimes (aI+b\hT^{m})] J$$
$$[aI\otimes I+b\hT^{m}\otimes \hT^{m} -(aI+b\hT^{m})\otimes (aI+b\hT^{m})]J=0$$
В спектральном представлении   действие $T^{(1,0)}$ заменяется на умножение на переменную $u$ (а $T^{(0,1)}$ -- на $w$). Получим
$$[(a-a^2)+(b-b^2){u_1u_2}  -ab(u_1+u_2)]Jf=0.$$
Отсюда находим однозначное  решение $u_1=\phi(u_2)$   уравнения 
$[(a-a^2)+(b-b^2){u_1u_2}  -ab(u_1+u_2)]=0$ (предположим, что оно существует). 
 Аналогично    $w_1=\phi(w_2)$. Следовательно, 
$\|Jf\|=\|\chi_\Gamma Jf\|$, где $\Gamma$ --  график аналитического отображения $\phi\otimes\phi$. Так как 
$\sigma\times\sigma (\Gamma)=0$ (спектральная мера действия  непрерывна), получаем $Jf=0$, $J=0$.

Таким образом,  наше исходное действие обладает однородным спектром кратности 2. Теорема 1 доказана.
\subsection{Набор кратностей $\{2,4\}$}

Рассмотрим  случай, когда  $\Z^2$-действие имеет набор кратностей $\{2,4\}$.

\begin{theorem}\label{theorem_n2}
    Для того, чтобы вполне эргодическое $\Z^2$-действие $\{R^z\otimes R^z\otimes T^z\otimes T^z\}_{z\in\Z^2}$
    обладало набором спектральных кратностей $\{2, 4\}$, достаточно выполнение следующих условий:
    \begin{enumerate}
        \item действия  $\{\hR^z\odot\hR^z\}_{z\in\Z^2}$ и $\{\hT^z\odot\hT^z\}_{z\in\Z^2}$  имеют простой спектр;
        \item для некоторых последовательностей $\{w_i\}$, $\{u_i\}$ $(i\to\infty)$ и 
     некоторых чисел $a_j, b_j>0$ и $j\in\{1,2\}$выполняется:
            \begin{equation} \label{t_condition}
                \begin{cases}
                    \hR^{w_i}\to a_1I+(1-a_1)\Theta, \\
                    \hT^{w_i}\to b_1\hT^{(1,0)}+(1-b_1)\Theta,
                \end{cases}
            \end{equation}
            \begin{equation} \label{r_condition}
                \begin{cases}
                    \hR^{u_i}\to a_2I+(1-a_2)\Theta, \\
                    \hT^{u_i}\to b_2\hT^{(0,1)}+(1-b_2)\Theta.
                \end{cases}
            \end{equation}
       
    \end{enumerate}
\end{theorem}

\proof{}
Пусть  $ U^z$   обозначает ограничение действия           $ \hR^z\odot \hR^z$         на  пространство, ортогональное константам,
а $V_z$ обозначает $ \hT^z\otimes \hT^z$.
Действие  $\{\hR^z\otimes \hR^z\otimes \hT^z\otimes \hT^z\}_{z\in\Z^2}$ изоморфно  
$$(\1+{\bf 2}U^z)\otimes  (\1+{\bf 2}V^z)= \1 + {\bf 2}U^z+{\bf 2}V^z+{\bf 4}U^z\otimes  V^z,$$
где $+$ обозначает  прямую сумму операторов,  ${\1}$ -  тождественный оператор на одномерном пространстве (констант),
${\bf 2}U^z:=U^z\oplus U^z$ и т.д. (подобные обозначения использовались в    \cite{Ryzhikov_MatSbornik_2009}).

Нам надо доказать, что  действия $U^z$, $V^z$, $U^z\otimes  V^z$  попарно спектрально дизъюнктны (унитарные представления дизъюнктны: нет ненулевого сплетения) и имеют   
простой спектр.  Это  автоматически приводит к очевидному набору кратностей $\{2,4\}$.
\\ 
Из условий теоремы имеем:

$U^{w_i}\to a_1^2I$,   \ $V^{w_i}\to b_1^2\hT^{(1,0)}\odot \hT^{(1,0)}$;\ \ \ 
$U^{u_i}\to a_2^2I$, \ $V^{u_i}\to b_2^2\hT^{(0,1)}\odot \hT^{(0,1)}$.
\\
Пусть $F$ -- циклический вектор для $U^{z}$, а   $G$ -- для $V^{z}$.  Тогда     $F\otimes G$ -- циклический вектор для  $U^z\otimes  V^z$.
Действительно, применяя шаг за шагом наши пределы,  получим, что векторы  

$F\otimes  \hT^{(1,0)}\odot \hT^{(1,0)} G$, \ \ $F\otimes  \hT^{(1,1)}\odot \hT^{(1,1) }G$, \ \ $F\otimes  \hT^{(m,n)}\odot \hT^{(m,n)} G$
\\
лежат в циклическом пространстве $C_{F\otimes G}$.  Следовательно, $U^z\otimes  V^z$ имеет простой спектр.

Дизъюнктность действий $U^z$ и $V^z$ вытекает из простоты спектра $U^z\otimes  V^z$  (иначе у $U^z\otimes  V^z$ появилась бы четная кратность). Дизъюнктность действий  $U^z$  и  $U^z\otimes  V^z$ доказывается так.
Если  
$$U^zJ=J(U^z\otimes  V^z),$$
то
$$U^{w_i}J=J(U^{w_i}\otimes  V^{w_i}),$$
$$a_1^2J=J(a_1^2I\otimes   b_1^2(\hT^{(1,0)}\odot \hT^{(1,0)})).$$
Так как $0<a_1, b_1<1,$  получим
$\|J\|=\|J\|b_1^2, \ \ J=0.$
Теорема 2 доказана.

\subsection{Набор кратностей $\{2, 4,\dots,2^k\}$}

\begin{lemma}
    Для того, чтобы действие $(\1+{\bf 2}V_1)\otimes \ldots\otimes  (\1+{\bf 2}V_k)$ имело 
 в качестве набора спектральных кратностей  множество $\{2, 4, 8, \ldots,2^k\}$,  достаточно, чтобы действие $(\1+V_1)\otimes \ldots\otimes  (\1+V_k)$ имело простой спектр.

\end{lemma}

Действительно, если прямая сумма 
$$ \1 { +} {\bf 1}(V_1+\ldots+V_k)+{\bf 1}(V_1\otimes V_2+\ldots + V_{k-1}\otimes V_k)+\\
\ldots + {\bf 1}V_1\otimes\ldots\otimes V_k$$
имеет простой спектр,  то  все слагаемые спектрально дизъюнктны, следовательно,
набор спектральных кратностей для действия 
$$ \1 + {\bf 2}(V_1+\ldots+V_k)+{\bf 4}(V_1\otimes V_2+\ldots + V_{k-1}\otimes V_k)+\\
\ldots + {\bf 2^k}V_1\otimes\ldots\otimes V_k$$
есть  множество $\{2, 4, 8, \ldots,2^k\}$.

\begin{theorem}\label{theorem_nk}. 
Вполне эргодическое действие 
$\{(\hT_1^z\otimes \hT_1^z)\otimes\ldots\otimes (\hT_k^z\otimes \hT_k^z)\}_{z\in\Z^2}$
имеет  набор спектральных кратностей   $\{2, 4, 8, \ldots,2^k\}$, если выполнены 
 следующие условия:
    \begin{enumerate}
        \item действия   $\hT_1^z\odot \hT_1^z,\ldots,\hT_k^z\odot \hT_k^z$ имеют простой спектр,
        \item для каждого $i\in\{1,2,\ldots,k\}$ существуют последовательности $\{w_n^i\}, \{u_n^i\}$ ($n\to\infty$) такие, что имеют место следующие слабые сходимости:
            \begin{equation} \label{t1_condition}
                \begin{cases}
                    \hT_1^{w_n^i}\to a_1^iI+(1-a_1^i)\Theta,\\
                                       \hT_2^{w_n^i}\to a_2^iI+(1-a_2^i)\Theta,\\
                    \ldots\\
                    \hT_i^{w_n^i}\to a_i^i\hT_i^{(1,0)}+(1-a_i^i)\Theta,\\
                    \ldots\\
                    \hT_k^{w_n^i}\to a_k^iI + (1-a_k^i)\Theta
                \end{cases}
                \begin{cases}
                    \hT_1^{u_n^i}\to b_1^iI+(1-b_1^i)\Theta,\\
                    \hT_2^{w_n^i}\to a_2^iI+(1-a_2^i)\Theta,\\
                    \ldots\\                    
\hT_i^{u_n^i}\to b_i^i\hT_i^{(0,1)}+(1-b_i^i)\Theta,\\
                    \ldots\\
                    \hT_k^{u_n^i}\to b_k^iI + (1-b_k^i)\Theta
                \end{cases}
            \end{equation}
            для некоторых чисел $a^i_j, b^i_j > 0$ и $j=\{1,\ldots,k\}$.
    \end{enumerate}
\end{theorem}

Доказательство. Действие $\{\hT_1^z\otimes \hT_1^z\otimes \hT_2^z\otimes \hT_2^z\otimes\ldots\otimes \hT_k^z\otimes \hT_k^z\}_{z\in\Z^2}$
    имеет набор кратностей $\{2, 4, 8, \ldots,2^k\}$,  если  действие 
$\{(\hT_1^z\odot \hT_1^z)\otimes (\hT_2^z\odot \hT_2^z)\otimes\ldots\otimes (\hT_k^z\odot \hT_k^z)\}_{z\in\Z^2}$ имеет простой спектр.  Действительно,  
действие  $\hT_i\otimes \hT_i$ в пространстве, ортогональном константам,   изоморфно прямой сумме двух копий $V_i:=\hT_i\odot \hT_i$ (что мы показали в параграфе 2), следовательно,  мы можем применить лемму 4.

  Слабые пределы обеспечивают и простоту спектра тензорного произведения 
$(\hT_1^z\odot \hT_1^z)\otimes\ldots\otimes (\hT_k^z\odot \hT_k^z)$  и на пространстве  $(L_2^0(\mu)\odot L_2^0(\mu))^{\otimes k}$  и  \\ дизъюнктность различных тензорных
произведений, входящих в него сомножителей.  Простота спектра доказывается  по индукции  методом    предыдущего параграфа.  Для этого в качестве $U$ рассматривается $(\hT_1^z\odot \hT_1^z)\otimes\ldots\otimes (\hT_{k-1}^z\odot \hT_{k-1}^z)$,  
а в качестве $V$  --  произведение $\hT_1^z\odot \hT_1^z$. 
 
Слабые пределы подобраны так, чтобы было легко показать дизъюнктность произведений различных наборов тензорных квадратов. Докажем, например,\\
дизъюнктность произведений $(\hT_1^z\otimes \hT_1^z)\otimes (\hT_2^z\otimes \hT_2^z)$  и
 $(\hT_2^z\otimes \hT_2^z)\otimes (\hT_3^z\otimes \hT_3^z)$.
                 Если $J$ сплетает  произведения, т.е.
 $$[(\hT_2^z\otimes \hT_2^z)\otimes (\hT_3^z\otimes \hT_3^z)]J=J
[(\hT_1^z\otimes \hT_1^z)\otimes (\hT_2^z\otimes \hT_2^z)],$$  то, 
учитывая  слабые  сходимости
                 $T_1^{w_n^3}\to a_1^3 I,$\ \
$T_2^{w_n^3}\to a_2^3 I,$\ \
                                      $ T_3^{w_n^3}\to a_3^3 T_3^{(1,0)}$, 
переходя к  пределам, получим 
$$  (a_2^3)^2  (a_3^3)^2 [(I\otimes I)\otimes (T_3^{(1,0)}\otimes T_3^{(1,0)})]J =  (a_1^3)^2(a_2^3)^2 J.$$  Мы оказались в ситуации, когда $(I\otimes V) J  = a J $, где оператор $V$  имеет непрерывный спектр,  следовательно, оператор $(I\otimes V)$ также не имеет собственных векторов, тем самым,  $J=0$.  Остальные случаи рассматриваются аналогично.

\end{document}